\documentclass[a4paper,11pt,twoside,reqno]{amsart}

\usepackage[utf8]{inputenc}
\usepackage[plainpages=false,pdfpagelabels=true]{hyperref}
\usepackage{amssymb,amsthm}
\usepackage[margin=1in]{geometry}

\numberwithin{equation}{section}

\newtheorem{Satz}{Theorem}[section]
\newtheorem{Prop}[Satz]{Proposition}
\newtheorem{Lem}[Satz]{Lemma}

\newtheorem{Cor}[Satz]{Corollary}
\theoremstyle{definition}
\newtheorem{Dfn}[Satz]{Definition}
\newtheorem{Bem}[Satz]{Remark}
\newtheorem{Bsp}[Satz]{Example}
\newcommand{\zarg}{d\phi(e_1)\wedge d\phi(e_2)}
\newcommand{\zargt}{d\phi_t(e_1)\wedge d\phi_t(e_2)}
\newcommand{\uarg}{du(e_1)\wedge du(e_2)}

\newcommand{\hol}{\operatorname{hol}}
\newcommand{\Hom}{\operatorname{Hom}}
\newcommand{\hess}{\operatorname{Hess}}
\newcommand{\cG}{{\mathcal G}}
\newcommand{\cL}{{\mathcal{L}}}
\newcommand{\tr}{\operatorname{Tr}}
\newcommand{\sff}{\mathrm{I\!I}}
\newcommand{\N}{\ensuremath{\mathbb{N}}}
\newcommand{\R}{\ensuremath{\mathbb{R}}}

\allowdisplaybreaks[1]

\renewcommand{\epsilon}{\varepsilon}

\title{The heat flow for the full bosonic string}
\author{Volker Branding}
\date{\today}
\address{TU Wien\\
Institut für diskrete Mathematik und Geometrie\\
Wiedner Hauptstraße 8–10, A-1040 Wien}
\email[]{volker@geometrie.tuwien.ac.at}

\begin{document}

\begin{abstract}
We study harmonic maps from surfaces coupled to a scalar and a two-form potential,
which arise as critical points of the action of the full bosonic string.
We investigate several analytic and geometric properties of these maps
and prove an existence result by the heat flow method. 
\end{abstract} 

\maketitle

\section{Introduction and Results}
Harmonic maps between Riemannian manifolds are one of the most studied variational problems
in differential geometry. Under the assumption that the target manifold has non-positive curvature Eells and Sampson 
established their famous existence result for harmonic maps making use of the heat flow method \cite{MR0164306}.
Moreover, in the case that the domain is two-dimensional harmonic maps belong to the class of conformally invariant variational problems 
yielding a rich structure.

However, harmonic maps from surfaces also have a dual life in theoretical physics.
More precisely, they arise as the \emph{Polyakov action} in bosonic string theory.
The full action for the bosonic string contains two additional terms, one of them being the pullback of a two-form from
the target and the other one being a scalar potential. It is the aim of this article to study the full action of the
bosonic string as a geometric variational problem.

There are already several mathematical results available for parts of the energy functional of
the full bosonic string:
In \cite{MR1433176} the authors consider the harmonic map energy together with a scalar potential.
The critical points of this energy functional are called \emph{harmonic maps with potential}.
One of the main results in that reference is that depending on the choice of potential,
the qualitative behavior of harmonic maps with potential differs from the one of harmonic maps.
There are several results available that characterize the properties of harmonic maps with potential:
This includes gradient estimates \cite{MR1667241} and Liouville theorems \cite{MR1618210}
for harmonic maps with potential from complete manifolds. In \cite{MR1680678} an existence result for harmonic maps
with potential from compact Riemannian manifolds with boundary is obtained, where it is assumed that the image of the map lies inside a convex ball.
An existence result for harmonic maps with potential to a target with negative curvature was obtained in \cite{MR1800592}
by the heat flow method. This result has been extended to the case of a domain manifold with boundary in \cite{MR1979036}.

Harmonic maps from surfaces coupled to a two-form potential have also been studied in the mathematical literature
since they give rise to the \emph{prescribed mean curvature equation}. 
Existence results via the heat flow for the prescribed mean curvature equation 
have been obtained for a flat target in \cite{MR1125012} and for a three-dimensional target with negative curvature in \cite{MR1424349}.

In theoretical physics the two-form potential is interpreted as giving rise to an external magnetic field.

Harmonic maps coupled to a two-form potential and spinor fields instead of a scalar potential have been studied in \cite{MR3305429}.

We will call the critical points of the full bosonic string action \emph{harmonic maps with scalar and two-form potential}
and we will generalize several results already obtained for harmonic maps and harmonic maps with potential.
Moreover, we will point out new phenomena that arise from the two-form potential.

This article is organized as follows. In Section 2 we analyze the energy functional of the full bosonic string and derive its critical points.
In Section 3 we study analytic and geometric aspects of the critical points.
In the last section we derive an existence result via the heat-flow method for both compact and non-compact target manifolds.

\section{The full bosonic string action}
Throughout this article \((M,h)\) is a Riemannian surface without boundary, we will mostly assume that \(M\) is compact,
and \((N,g)\) a closed, oriented Riemannian manifold of dimension \(\dim N\geq 3\).
For a map \(\phi\colon M\to N\) we consider the square of its differential giving rise to the usual harmonic energy.
Let \(B\) be a two-form on \(N\), which we pull back by the map \(\phi\).
In addition, let \(V\colon N\to\R\) be a scalar function, which we mostly assume to be smooth. By \(R\) we denote the scalar curvature of the domain \(M\).

In the physical literature the full action for the bosonic string is given by
\begin{equation}
\label{energy-functional}
E(\phi)=\int_M\big(\frac{1}{2}|d\phi|^2+\phi^\ast B+RV(\phi)\big)dM,
\end{equation}
see for example \cite{MR2151029}, p. 108.
Due to the uniformization theorem we can assume that the scalar curvature \(R\) on the domain \(M\) is constant.
In string theory, the potential \(V(\phi)\) is usually referred to as \emph{dilaton field}.

For the sake of completeness we want to mention that in theoretical physics one also defines an action functional 
that locally has the form \eqref{energy-functional} but has a different global structure.
More precisely, one replaces the two-form contribution in \eqref{energy-functional} by an object that locally
looks like a two-form, but forms a more general object from the global point of view. 
One then studies the \(U(1)\)-valued functional
\begin{equation}
\label{energy-functional-u1}
\exp(iE(\phi)):=\exp\big(i\int_M(\frac{1}{2}|d\phi|^2+RV(\phi))dM\big)\cdot\hol(\phi),
\end{equation}
where \(\hol(\phi)\) denotes the holonomy of the gerbe \(\cG\) along the map \(\phi\).
By the three-form \(\Omega\) we denote the globally defined curvature associated to the gerbe \(\cG\).
For a mathematical introduction to gerbes and surface holonomy see \cite{MR2648325}.
This gives rise to the so-called \emph{B-field action} in string theory, which is defined as
\begin{equation*}
E_B(\phi):=-i\log(\hol(\phi)) 
\end{equation*}
such that \(\hol(\phi)=\exp(iE_B(\phi))\).

In the case that the gerbe \(\cG\) is trivial, which corresponds to the fact that the three-form \(\Omega\)
is exact, we can define the \emph{B-field action} as
\begin{equation*}
E_B(\phi):=\int_M\phi^\ast BdM
\end{equation*}
leading to the functional \eqref{energy-functional}. 
From an analytical point of view the \(U(1)\)-valued functional \eqref{energy-functional-u1} is more difficult,
and we will restrict ourselves to the functional \eqref{energy-functional}.

Let us derive the critical points of \eqref{energy-functional}.
\begin{Prop}
The Euler-Lagrange equation of the functional \eqref{energy-functional} is given by
\begin{equation}
\label{euler-lagrange}
\tau(\phi)=Z(\zarg)+R\nabla V(\phi),
\end{equation}
where \(\tau(\phi)\) denotes the tension field of the map \(\phi\) and the vector-bundle homomorphism \(Z \in \Gamma (\Hom(\Lambda^2T^\ast N,TN))\) is defined by the equation
\begin{equation}
\label{def-Z}
\Omega(\eta,\xi_1,\xi_2)=\langle Z(\xi_1\wedge\xi_2),\eta\rangle,
\end{equation}
where \(\Omega=dB\) is a three-form on \(N\) and \(e_1,e_2\) an orthonormal basis of \(TM\).

\end{Prop}
\begin{proof}
We consider a smooth variation of \(\phi\) that is \(\phi_t\colon (-\epsilon,\epsilon)\times M\to N\) with \(\frac{\partial\phi_t}{\partial t}\big|_{t=0}=\eta\).
It is well-known that
\[
\frac{d}{dt}\big|_{t=0}\frac{1}{2}\int_M|d\phi_t|^2dM=-\int_M\langle\tau(\phi),\eta\rangle dM, \qquad \frac{d}{dt}\big|_{t=0}\int_M RV(\phi_t)dM=\int_M R\langle\nabla V(\phi),\eta\rangle dM.
\]
To calculate the variation of the two-form \(B\) we choose an orthonormal basis \(e_1,e_2\) of \(TM\) and calculate
\begin{align*}
\frac{d}{dt}(\phi_t^\ast B)(e_1,e_2)=&\phi_t^\ast(\cL_{\partial_t}B)(e_1,e_2) \\
=&\phi_t^\ast\big(\iota_{\partial_t}dB(e_1,e_2)+d(\iota_{\partial_t} B(e_1,e_2))\big) \\
=&\phi_t^\ast\big(\Omega(\partial_t,e_1,e_2))+d(\phi_t^\ast(\iota_{\partial_t}B(e_1,e_2))),
\end{align*}
where we used that \(dB=\Omega\). Here, \(\cL\) denotes the Lie-derivative applied to differential forms.
Moreover, we applied several identities for derivatives of differential forms, which can be found in \cite{MR3012162}, pp. 170-174. 

Integrating over \(M\) we thus obtain
\[
\frac{d}{dt}\big|_{t=0}\int_M(\phi^\ast_tB)(e_1,e_2)dM=\int_M\Omega(\eta,d\phi(e_1),d\phi(e_2))dM.
\]
Using the vector-bundle homomorphism \eqref{def-Z}, we thus find
\[
\frac{d}{dt}\big|_{t=0}E(\phi_t)=\int_M\langle\eta,-\tau(\phi)+Z(\zarg)+R\nabla V(\phi)\rangle dM
\]
yielding the claim.
\end{proof}
We call solutions of \eqref{euler-lagrange} \emph{harmonic maps with scalar and two-form potential}.

\begin{Bem}
Note that the energy functionals \eqref{energy-functional} and \eqref{energy-functional-u1} have the same critical points since
\[
\frac{d}{dt}\big|_{t=0}\exp(iE(\phi_t))=0\Leftrightarrow \frac{d}{dt}\big|_{t=0}\int_M(\frac{1}{2}|d\phi_t|^2+RV(\phi_t))dM+\frac{d}{dt}\big|_{t=0}E_B(\phi_t)=0.
\]
\end{Bem}

Whenever using local coordinates,
we will use greek indices on the domain \(M\) and latin indices on the target \(N\). Moreover,
we will make use of the Einstein summation convention, that is, we sum over repeated indices.
In terms of local coordinates \(x_\alpha\) on \(M\) and \(y^i\) on \(N\) the Euler-Lagrange equation reads
\begin{equation}
\label{euler-lagrange-local}
\Delta\phi^i=-\Gamma^i_{jk}\frac{\partial\phi^j}{\partial x_\alpha}\frac{\partial\phi^k}{\partial x_\beta}h_{\alpha\beta}
+Z^i(\partial_{y^j}\wedge\partial_{y^k})\frac{\partial\phi^j}{\partial x_1}\frac{\partial\phi^k}{\partial x_2}
+Rg^{ik}\frac{\partial V}{\partial y^k}.
\end{equation}

\begin{Bem}
The vector-bundle homomorphism \(Z\) can also be interpreted as arising from a metric connection with totally antisymmetric torsion.
In this case one has
\[
\nabla^{Tor}_XY=\nabla^{LC}_XY+A(X,Y),
\]
where \(\nabla^{LC}\) denotes the Levi-Cevita connection, \(X,Y\) are vector fields and \(A(X,Y)\) is a skew-adjoint endomorphism.
The endomorphism \(A(X,Y)\) satisfies
\[
\langle A(X,Y),Z\rangle=\Omega(X,Y,Z)
\]
with \(\Omega\in\Gamma(\Lambda^3T^\ast N)\) similar to \eqref{def-Z}. For more details on metric connections with torsion
see \cite{MR3493217} and references therein.
\end{Bem}

\begin{Bem}
In principle one could also study the functional \eqref{energy-functional} for a higher-dimensional domain \(M\) with \(m=\dim M\geq 2\).
However, then one needs to pull back an \(m-\)form from the target leading to an Euler-Lagrange equation with a higher nonlinearity
on the right hand side, see \cite{KohDiss}, Chapter 2, for a detailed analysis.
\end{Bem}

\begin{Lem}[Second Variation]
Let \(\phi\colon M\to N\) be a harmonic map with scalar and two-form potential.
The second variation of the energy functional \eqref{energy-functional} is given by
\begin{align}
\label{second-variation}
\frac{\partial^2}{\partial t^2}\big|_{t=0}E(\phi_t)
=\int_M&\big(|\nabla\xi|^2-\langle R^N(d\phi,\xi)d\phi,\xi\rangle+\langle\xi,(\nabla_\xi Z)(\zarg)\rangle \\
\nonumber+&\langle\xi,Z(\nabla\xi(e_1)\wedge d\phi(e_2))\rangle+\langle\xi,Z(d\phi(e_1)\wedge \nabla\xi(e_2))\rangle
+R\hess V(\xi,\xi)\big)dM,
\end{align}
where \(\xi=\frac{\partial\phi_t}{\partial t}\big|_{t=0}\).
\end{Lem}
\begin{proof}
This follows by a direct calculation.
\end{proof}

\begin{Bsp}
Suppose that \(\phi\) is an isometric immersion from a closed oriented Riemann surface \((M,h)\) into an oriented Riemannian three manifold \((N,g)\).
In this case the tension field \(\phi\) is related to the mean curvature vector \(H(\phi)\) by
\[
\frac{1}{2}\tau(\phi)=H(\phi).
\]
Any three-form \(\Omega\) on \(N\) must be a multiple of the volume form \(\operatorname{vol}_g\), that is
\[
\Omega=f\operatorname{vol}_g
\]
for some smooth function \(f\colon N\to\R\).
Thus, if we denote the unit normal field of \(\phi(M)\subset N\) by \(\nu\) the Euler-Lagrange equation \eqref{euler-lagrange}
is equivalent to 
\begin{equation}
\label{example-three-mfd}
2H(\phi)=\pm (f\nu)\circ\phi+R\nabla V(\phi).
\end{equation}
The sign in front of the first term on the right hand side depends on whether the map \(\phi\) is
orientation preserving or not, see \cite{KohDiss}, Remark 2.7.

In case that the function \(f\) is constant and \(\nabla V\sim\nu\) the equation \eqref{example-three-mfd} has some similarity
with the equation for \emph{linear Weingarten surfaces}. These are surfaces immersed in a three-dimensional manifold satisfying
\[
aH+bK=c,
\]
where \(H\) denotes the mean curvature, \(K\) the Gauss curvature and \(a,b,c\) are non-zero real numbers.
\end{Bsp}

\section{Properties of Harmonic maps with scalar and two-form potential}
In this section we study several properties of solutions of \eqref{euler-lagrange}.
\begin{Lem}
Let \(\phi\in C^2(M,N)\) be a solution of \eqref{euler-lagrange}. Then the following formulas hold:
\begin{align}
\Delta\frac{1}{2}|d\phi|^2=&|\nabla d\phi|^2-\langle R^N(d\phi(e_\alpha),d\phi(e_\beta))d\phi(e_\alpha),d\phi(e_\beta)\rangle
+\langle d\phi(\operatorname{Ric}^M(e_\alpha)),d\phi(e_\alpha)\rangle \\
\nonumber &-\langle Z(\zarg),\tau(\phi)\rangle+R\hess V(d\phi,d\phi) \\
\intertext{and}
\label{bochner-vphi}
\Delta (V\circ\phi)&=dV\big(R\nabla V(\phi)+Z(\zarg)\big)+\tr\hess V(d\phi,d\phi),
\end{align}
where \(\operatorname{Ric}^M\) denotes the Ricci-curvature on \(M\) and \(R^N\) the curvature tensor on \(N\).
\end{Lem}

\begin{proof}
This follows by a direct calculation.
\end{proof}

\begin{Dfn}
We define a two-tensor by 
\begin{equation}
\label{stress-energy tensor}
S_{\alpha\beta}:=\frac{1}{2}h_{\alpha\beta}|d\phi|^2-\langle d\phi(e_\alpha),d\phi(e_\beta)\rangle+RV(\phi)h_{\alpha\beta}.
\end{equation}
The tensor \(S_{\alpha\beta}\) is called the \emph{stress-energy tensor}.
\end{Dfn}
\begin{Lem}
The stress-energy tensor \(S_{\alpha\beta}\) is divergence free and symmetric if \(\phi\)
is a harmonic map with scalar and two-form potential.
\end{Lem}
\begin{proof}
By a direct calculation we find
\begin{align*}
\nabla_{e_\alpha}S_{\alpha\beta}=&-\langle\tau(\phi),d\phi(e_\beta)\rangle+\langle R\nabla V(\phi),d\phi(e_\beta)\rangle \\
&=-\langle Z(\zarg),d\phi(e_\beta)\rangle \\
&=-\Omega(d\phi(e_1),d\phi(e_2),d\phi(e_\beta))=0,
\end{align*}
which proves the claim.
\end{proof}
Note that \(S_{\alpha\beta}\) is no longer trace-free, which corresponds to the fact that the scalar potential \(V(\phi)\)
in the action functional does not respect the conformal symmetry.

We can express the stress-energy tensor \eqref{stress-energy tensor} invariantly as
\begin{align}
\label{stress-energy tensor-invariant}
S_V(\phi)=e(\phi)h-\phi^\ast g+RV(\phi)h
\end{align}
with the energy density \(e(\phi)=\frac{1}{2}|d\phi|^2\).

With the help of \eqref{stress-energy tensor-invariant} we now derive a monotonicity formula for \(e(\phi)+RV(\phi)\).
A similar calculation for harmonic maps with scalar potential has been carried out in \cite{MR2929724}.
For the monotonicity formula for harmonic maps we refer to the book \cite{MR1391729}, for a general treatment of stress-energy tensors
with applications to harmonic maps we refer to \cite{MR2453667}.
Note that we do not have to assume that \(M\) is compact to derive the monotonicity formula.

Let \((M,h_0)\) be a complete Riemannian surface with a pole \(x_0\). Let \(r(x)\) be the Riemannian distance
function relative to the point \(x_0\). By \(\lambda_i,i=1,2\) we denote the eigenvalues of \(\hess_{h_0}(r^2)\).

\begin{Prop}[Monotonicity formula]
Let \(\phi\colon (M,f^2h_0)\to (N,g)\) be a harmonic map with scalar and two-form potential, where \(M\) is a complete Riemannian surface and \(N\) a Riemannian manifold.
Suppose that 
\begin{equation}
\label{assumptions-monotonicity}
RV(\phi)> 0, \qquad \qquad r\frac{\partial\log f}{\partial r}\geq 0, \qquad \qquad \frac{1}{2}\big(\sum_{i=1}^2\lambda_i-2\lambda_{max}\big)\geq\sigma
\end{equation}
for a positive constant \(\sigma\).
Then the following inequality holds
\begin{equation}
\frac{\int_{B_{\rho_1}(x_0)}\big(e(\phi)+RV(\phi)\big)}{\rho_1^\sigma}\leq \frac{\int_{B_{\rho_2}(x_0)}\big(e(\phi)+RV(\phi)\big)}{\rho_2^\sigma}
\end{equation}
for any \(0<\rho_1\leq\rho_2\). Here, \(B_\rho(x)\) denotes the geodesic ball with radius \(\rho\) around the point \(x\).
\end{Prop}

The proof is similar to the proof of Theorem 4.1 in \cite{MR2929724}.

\begin{proof}
For a symmetric \((2,0)\)-tensor \(S\) the following identity holds
\begin{align}
\label{tensor-identity}
\operatorname{div}(\iota_X S)=(\operatorname{div}S)(X)+\frac{1}{2}\langle S,\cL_Xh\rangle.
\end{align}
Here, \(\cL_X\) represents the Lie derivative with respect to the vector field \(X\).
Integrating the divergence of \eqref{stress-energy tensor-invariant} over the ball \(B_\rho(x_0)\) with radius \(\rho\) around the point \(x_0\), using \eqref{tensor-identity} and the fact
that the stress-energy-tensor is divergence free, we obtain
\begin{align*}
\int_{\partial B_{\rho}(x_0)}S_V(\phi)(X,\nu)=\frac{1}{2}\int_{B_{\rho}(x_0)}\langle S_V(\phi),\cL_Xh\rangle.
\end{align*}
We choose \(X=r\frac{\partial}{\partial r}\) and get
\begin{align*}
\frac{1}{2}\cL_Xh=\frac{1}{2}\cL_X(f^2h_0)=rf\frac{\partial f}{\partial r}h_0+\frac{1}{2}f^2\cL_Xh_0,
\end{align*}
where we used that the metric on \(M\) is given by \(h=f^2h_0\). Hence, we find
\begin{align*}
\frac{1}{2}\langle S_V(\phi),\cL_Xh\rangle=2RV(\phi)r\frac{\partial\log f}{\partial r}+\frac{1}{2}f^2\langle S_V(\phi),\hess_{h_0}(r^2)\rangle
\end{align*}
since \(\dim M=2\).

To estimate the last term on the right hand side, we choose a local basis \(e_1,e_2\) on \(M\) with respect to \(h_0\) such that \(\hess_{h_0}(r^2)\)
becomes a diagonal matrix with respect to \(e_1,e_2\). Then \(\tilde{e}_\alpha:=f^{-1}e_\alpha,\alpha=1,2\) is an orthonormal basis with respect to \(h\).
Thus, we obtain
\begin{align*}
f^2\langle S_V(\phi),\hess_{h_0}(r^2)\rangle=&f^2\sum_{\alpha,\beta=1}^2S_V(\phi)(\tilde{e}_\alpha,\tilde{e}_\beta)\hess_{h_0}(r^2)(\tilde{e}_\alpha,\tilde{e}_\beta) \\
=&\sum_{\alpha,\beta=1}^2\big(e(\phi)h(\tilde{e}_\alpha,\tilde{e}_\beta)-\phi^\ast g(\tilde{e}_\alpha,\tilde{e}_\beta)+RV(\phi)g(\tilde{e}_\alpha,\tilde{e}_\beta)\big)\hess_{h_0}(r^2)(e_\alpha,e_\beta) \\
=&\sum_{\alpha=1}^2(e(\phi)+RV(\phi))\hess_{h_0}(r^2)(e_\alpha,e_\alpha)
-\sum_{\alpha=1}^2\phi^\ast g(\tilde{e_\alpha},\tilde{e_\alpha})\hess_{h_0}(r^2)(e_\alpha,e_\alpha) \\
\geq &\big(e(\phi)+RV(\phi)\big)\sum_{i=1}^2\lambda_i-2e(\phi)\lambda_{max} \\
\geq &\big(e(\phi)+RV(\phi)\big)\big(\Lambda -2\lambda_{max}\big),
\end{align*}
where \(\Lambda=\lambda_1+\lambda_2\).
Hence, we obtain
\begin{align*}
\frac{1}{2}\langle S_V(\phi),\cL_xh\rangle\geq 2RV(\phi)r\frac{\partial\log f}{\partial r}+\frac{1}{2}\big(e(\phi)+RV(\phi)\big)\big(\Lambda -2\lambda_{max}\big).
\end{align*}
On the other hand, by the coarea formula and \(|\nabla r|=f^{-1}\) we find
\begin{align*}
\nonumber\int_{\partial B_{\rho}(x_0)}S_V(\phi)(r\frac{\partial}{\partial r},\nu)\leq &\int_{\partial B_{\rho}(x_0)}(e(\phi)+RV(\phi))h(r\frac{\partial}{\partial r},\nu)\\
\nonumber=&\rho\int_{\partial B_{\rho}(x_0)}(e(\phi)+RV(\phi))f\\
\nonumber=&\rho\frac{d}{d\rho}\int_{B_\rho(x_0)}\big(\int_{\partial B_{r}(x_0)}\frac{e(\phi)+RV(\phi)}{|\nabla r|}\big)dr \\
=&\rho\frac{d}{d\rho}\int_{B_\rho(x_0)}\big(e(\phi)+RV(\phi)\big).
\end{align*}
Thus, by the assumptions \eqref{assumptions-monotonicity} we obtain the following inequality
\begin{align*}
\rho\frac{d}{d\rho}\int_{B_\rho(x_0)}(e(\phi)+RV(\phi))\geq\int_{B_\rho(x_0)}\frac{1}{2}\big(e(\phi)+RV(\phi)\big)\big(\Lambda -2\lambda_{max}\big).
\end{align*}
Again, by assumption we get
\begin{equation*}
\rho\frac{d}{d\rho}\int_{B_\rho(x_0)}(e(\phi)+RV(\phi))\geq\sigma\int_{B_\rho(x_0)}(e(\phi)+RV(\phi)).
\end{equation*}
This can be rewritten as
\[
\frac{d}{d\rho}\frac{\int_{B_\rho(x_0)}(e(\phi)+RV(\phi))}{\rho^\sigma}\geq 0
\]
and the claim follows by integration with respect to \(\rho\).
\end{proof}

As for harmonic maps (\cite{MR510549}, Theorem 2) and harmonic maps with potential (\cite{MR1433176}, Proposition 2)
we can prove a unique continuation theorem for harmonic maps with scalar and two-form potential.
To obtain this result we recall the following (\cite{MR0092067}, p. 248)
\begin{Satz}
\label{aro-theorem}
Let \(A\) be a linear elliptic second-order differential operator defined on a domain \(D\) of \(\R^n\).
Let \(u=(u^1,\ldots,u^n)\) be functions in \(D\) satisfying the inequality
\begin{equation}
\label{aro-voraus}
|Au^j|\leq C\big(\sum_{\alpha,i}\big|\frac{\partial u^i}{\partial x^\alpha}\big|+\sum_i|u^i|\big).
\end{equation}
If \(u=0\) in an open set, then \(u=0\) throughout \(D\).
\end{Satz}

Making use of this result we can prove the following
\begin{Prop}
Let \(\phi,\phi'\in C^2(M,N)\) be two harmonic maps with scalar and two-form potential. Moreover, assume that \(V\colon N\to\R\)
is a \(C^{1,1}\) function. If \(\phi\) and \(\phi'\) are equal on a connected open set \(W\) of \(M\)
then they coincide on the whole connected component of \(M\) which contains \(W\).
\end{Prop}
\begin{proof}
Let \(U\) be a coordinate ball on \(M\) such that \(\phi=\phi'\) in some open subset.
By shrinking \(U\) if necessary we can assume that both \(\phi\) and \(\phi'\) map \(U\)
into a single coordinate chart in \(N\). We write \(y^i(x)\) for \(y^i(\phi(x))\) and \(z^i(x)\) for \(z^i(\phi'(x))\).
We consider the function \(u^i:=y^i-z^i\). Using the local form of the Euler-Lagrange equation \eqref{euler-lagrange-local} we find
\begin{align*}
\Delta u^i=&-\Gamma^i_{jk}(y)\frac{\partial y^j}{\partial x_\alpha}\frac{\partial y^k}{\partial x_\beta}h_{\alpha\beta}
+\Gamma^i_{jk}(z)\frac{\partial z^j}{\partial x_\alpha}\frac{\partial z^k}{\partial x_\beta}h_{\alpha\beta} \\
&-Z^i(y)(\partial_{y^j}\wedge\partial_{y^k})\frac{\partial y^j}{\partial x_1}\frac{\partial y^k}{\partial x_2}
+Z^i(z)(\partial_{y^j}\wedge\partial_{y^k})\frac{\partial z^j}{\partial x_1}\frac{\partial z^k}{\partial x_2} \\
&+Rg^{ik}(y)\frac{\partial V}{\partial y^k}(y)-Rg^{ik}(z)\frac{\partial V}{\partial y^k}(z).
\end{align*}
By adding several zero's on the right hand side the above equation acquires the form \eqref{aro-voraus} and
the result follows by application of Theorem \ref{aro-theorem}. For more details see the proof of Theorem 2
in \cite{MR510549}.
\end{proof}

\begin{Bem}
In the case that the scalar potential \(V(\phi)\) vanishes identically, the energy functional \eqref{energy-functional}
is conformally invariant. In this case one can exploit the conformal invariance to prove that critical points
cannot have isolated singularities, whenever a certain energy is finite. 
This follows from the main theorem in \cite{MR744314}, where a removable singularity theorem for critical points
of a large class of conformally invariant energy functionals is established.
\end{Bem}

As a next step we use the embedding theorem of Nash to isometrically embed \(N\) into some \(\R^q\).
We denote the isometric embedding by \(\iota\) and consider the composite map \(u:=\iota\circ\phi\colon M\to\R^q\).
By \(\tilde{N}\) we denote the tubular neighborhood of \(\iota(N)\subset\R^q\).
Let \(\pi\colon\tilde{N}\to\iota(N)\) be the canonical projection  which assigns to each \(z\in\tilde{N}\) the closest point
in \(\iota(N)\) from \(z\).

\begin{Lem}
\label{embedding-rq}
Assume that \(N\subset\R^q\). Then the Euler-Lagrange equation \eqref{euler-lagrange} acquires the form
\begin{equation}
\Delta u=\sff(du,du)+\tilde{Z}(du(e_1)\wedge du(e_2))+R\widetilde{\nabla V}(u),
\end{equation}
where \(\sff(du,du):=\tr\nabla d\pi(du,du)\). Moreover, \(\tilde{Z}\) and \(\widetilde{\nabla V}\) 
denote the extensions of \(Z\) and \(\nabla V\) to the ambient space \(\R^q\).
\end{Lem}
\begin{proof}
This follows from the chain rule for the tension field of composite maps, that is 
\[
\Delta (\iota\circ\phi)=d\iota(\tau(\phi))+\tr\nabla d\pi(d\phi,d\phi).
\]
The homomorphism \(Z\) and \(\nabla V\) can be extended to the ambient space by projecting
to a tubular neighborhood, for more details see \cite{MR2551140}, p. 463 and \cite{MR1800592}, p. 557.
\end{proof}

\section{Existence Results via the heat flow}
In this section we derive an existence result for critical points of \eqref{energy-functional}
by the heat flow method. In order to achieve this result we will assume that \(\Omega\) is exact
such that we have a variational structure that enables us to derive the necessary estimates
for convergence of the gradient flow. A similar approach for geodesics coupled to a magnetic field 
was performed in \cite{magnetic-geodesics}.

In order to control the non-linearities arising from the two-form 
we will have to restrict to target spaces with negative sectional curvature. 
A similar idea has been used in \cite{MR1424349} for the heat flow of the prescribed mean curvature equation.
More precisely, we will use the negative curvature of the target to control the non-linearities arising from
the two-form potential.

For a general introduction to harmonic maps and their heat flows see the book \cite{MR2431658}.

The gradient flow of the functional \eqref{energy-functional} is given by
\begin{equation}
\label{flow-mfd}
\frac{\partial\phi_t}{\partial t}(x,t)=\big(\tau(\phi_t)-Z(\zargt)-R\nabla V(\phi_t)\big)(x,t),
\end{equation}
where \(\phi_t\colon M\times [0,T)\to N\) and initial data \(\phi(x,0)=\phi_0(x)\). 
To establish the short-time existence of \eqref{flow-mfd}
we use the Nash embedding theorem to isometrically embed \(N\) into \(\R^q\).

\begin{Lem}
Suppose that \(N\) is isometrically embedded into \(\R^q\). Then \eqref{flow-mfd} acquires the form
\begin{equation}
\label{flow-rq}
\frac{\partial u}{\partial t}=\Delta u-\sff(du,du)-Z(\uarg)-R\nabla V(u),
\end{equation} 
where \(u\colon M\times [0,T)\to\R^q\).
\end{Lem}
\begin{proof}
This follows from Lemma \ref{embedding-rq}, we will omit the tildes in order not to blow up the notation.
\end{proof}
As a first step, we establish the existence of a short-time solution.
\begin{Lem}
For \(\phi_0\in C^{2+\alpha}(M,N)\) and \(V\in C^{1,1}(N)\) there exists a unique, smooth solution to \eqref{flow-mfd} for \(t\in [0,T_{max})\).
\end{Lem}
\begin{proof}
This can be proven using the Banach fixed point theorem, which requires that the potential \(V\in C^{1,1}(N)\).
For more details, see \cite{MR2431658}, Chapter 5.
\end{proof}

To extend the solution beyond \(T_{max}\) we will make use of the following Bochner formulae:

\begin{Lem}
Let \(\phi_t\colon M\times [0,T_{max})\to N\) be a smooth solution of \eqref{flow-mfd}. Then we have for all \(t\in [0,T_{max})\)

\begin{align}
\label{bochner-dphi}
  \frac{\partial}{\partial t}\frac{1}{2}|d\phi_t|^2=&\Delta \frac{1}{2}|d\phi_t|^2 -|\nabla d\phi_t|^2
 +\langle R^N(d\phi_t(e_\alpha),d\phi_t(e_\beta))d\phi_t(e_\alpha),d\phi_t(e_\beta)\rangle\\
\nonumber& -\langle d\phi_t(Ric^M(e_\alpha)),d\phi_t(e_\alpha)\rangle+\langle Z(\zargt),\tau(\phi_t)\rangle-R\hess V(d\phi_t,d\phi_t) \\
\intertext{and}
\label{bochner-phit}
\frac{\partial}{\partial t}\frac{1}{2}|\frac{\partial\phi_t}{\partial t}|^2=&
\Delta\frac{1}{2}|\frac{\partial\phi_t}{\partial t}|^2-|\nabla\frac{\partial\phi_t}{\partial t}|^2
+\langle R^N(d\phi_t(e_\alpha),\frac{\partial\phi_t}{\partial t})d\phi_t(e_\alpha),\frac{\partial\phi_t}{\partial t}\rangle\\
\nonumber&-\langle\frac{\nabla}{\partial t}Z(\zargt),\frac{\partial\phi_t}{\partial t}\rangle
-R\hess V(\frac{\partial\phi_t}{\partial t},\frac{\partial\phi_t}{\partial t}).
\end{align}
\end{Lem}
\begin{proof}
This follows from the standard Bochner formulas using 
\[
\langle\nabla Z(\zargt),d\phi_t\rangle=-\langle Z(\zargt),\tau(\phi_t)\rangle.
\]
\end{proof}
The Bochner formulae will be the key-tool to achieve long-time existence and convergence of the evolution equation \eqref{flow-mfd}.
However, we have to distinguish between the cases of a compact and a non-compact target.

\subsection{Compact Target}
For a compact target manifold \(N\) we will prove the following
\begin{Satz}
\label{theorem-compact}
Let \((M,h)\) be a closed Riemannian surface and \((N,g)\) be a closed, oriented Riemannian manifold
with negative sectional curvature. Moreover, suppose that \(\Omega\) is exact, \(|B|_{L^\infty}<1/2\), \(V\in C^{2,1}(N)\) and
that the homomorphism \(Z\) satisfies
\begin{equation*}
\frac{1}{2}|Z|^2_{L^\infty}\leq\kappa_N,
\end{equation*}
where \(\kappa_N\) denotes and upper bound on the sectional curvature on \(N\).
Then \(\phi_t\colon M\times [0,\infty)\to N\) 
subconverges as \(t\to\infty\) to a harmonic map with scalar and two-form potential \(\phi_\infty\),
which is homotopic to \(\phi_0\).
\end{Satz}
We will divide the proof of Theorem \ref{theorem-compact} into several steps.

\begin{Lem}
Let \(\phi_t\colon M\times [0,T_{max})\to N\) be a smooth solution of \eqref{flow-mfd} and \(V\in C^2(N)\). Then for all \(t\in [0,T_{max})\)
the following inequalities hold:
\begin{align}
\label{estimate-dphi}
  \frac{\partial}{\partial t}\frac{1}{2}|d\phi_t|^2\leq&\Delta \frac{1}{2}|d\phi_t|^2
  +c_1|d\phi_t|^2+(\frac{1}{2}|Z|_{L^\infty}^2-\kappa_N)|d\phi_t|^4
\intertext{and}
\label{estimate-phit}
\frac{\partial}{\partial t}\frac{1}{2}|\frac{\partial\phi_t}{\partial t}|^2\leq&
\Delta\frac{1}{2}|\frac{\partial\phi_t}{\partial t}|^2
+\big(|\nabla Z|_{L^\infty}+\frac{1}{4}|Z|_{L^\infty}^2-\kappa_N\big)|d\phi_t|^2|\frac{\partial\phi_t}{\partial t}|^2+c_2|\frac{\partial\phi_t}{\partial t}|^2
\end{align}
with \(c_1:=|\operatorname{Ric}|_{L^\infty}+|R|_{L^\infty}|\hess V|_{L^\infty}, c_2:=|R|_{L^\infty}|\hess V|_{L^\infty}\) and \(\kappa_N\) denotes an upper bound on
the sectional curvature on \(N\).
\end{Lem}

\begin{proof}
For the first assertion we consider \eqref{bochner-dphi} and estimate
\[
|\langle Z(\zargt),\tau(\phi_t)\rangle|\leq \sqrt{2}|Z|_{L^\infty}|d\phi_t|^2|\nabla d\phi_t|,
\]
which yields
\[
\langle Z(\zargt),\tau(\phi_t)\rangle-|\nabla d\phi_t|^2\leq\frac{1}{2}|Z|^2_{L^\infty}|d\phi_t|^4.
\]
By assumption \(N\) is compact and we can estimate the Hessian of the potential \(V(\phi)\) by its maximum yielding
the first claim.

For the the second claim we consider \eqref{bochner-phit} and calculate
\begin{align*}
\frac{\nabla}{\partial t}Z(\zargt)=&(\nabla_{d\phi_t(\partial_t)}Z)(\zargt)+Z(\frac{\nabla}{\partial t}d\phi_t(e_1)\wedge d\phi_t(e_2)) \\
&+Z(d\phi_t(e_1)\wedge\frac{\nabla}{\partial t}d\phi_t(e_2)),
\end{align*}
which gives
\begin{align*}
\langle\frac{\nabla}{\partial t}Z(\zargt),\frac{\partial\phi_t}{\partial t}\rangle
\leq&|\nabla Z|_{L^\infty}|d\phi_t|^2|\frac{\partial\phi_t}{\partial t}|^2+|Z|_{L^\infty}|\frac{\nabla}{\partial t}d\phi_t||d\phi_t||\frac{\partial\phi_t}{\partial t}|.
\end{align*}
The result follows by applying Young's inequality to the last term on the right hand side.
\end{proof}

Via the maximum principle we thus obtain the following
\begin{Cor}
Let \(\phi_t\colon M\times [0,T_{max})\to N\) be a smooth solution of \eqref{flow-mfd}.
If \(\frac{1}{2}|Z|_{L^\infty}^2\leq\kappa^N\) then for all \(t\in [0,T_{max})\) the following estimates hold:
\begin{align}
\label{bound-dphi} |d\phi_t|^2&\leq |d\phi_0|^2e^{2c_1t}, \\
\label{bound-phit} |\frac{\partial\phi_t}{\partial t}|^2&\leq |\frac{\partial\phi_0}{\partial t}|^2e^{\frac{|\nabla Z|_{L^\infty}|d\phi_0|^2}{c_1}e^{2c_1t}+c_2t}.
\end{align}
\end{Cor}
\begin{proof}
Using the assumptions the first statement follows by applying the maximum principle to \eqref{estimate-dphi}.
For the second statement we apply the maximum principle to \eqref{estimate-phit} using the bound \eqref{bound-dphi}.
\end{proof}

\begin{Lem}
\label{lemma-regularity}
Let \(u\colon M\times [0,T_{max})\to\R^q\) be a smooth solution of \eqref{flow-rq}.
If \(\frac{1}{2}|Z|_{L^\infty}^2\leq\kappa^N\) then the following inequality holds
\begin{align*}
|u(\cdot,t)|_{C^{2+\alpha}(M,N)}+\big|\frac{\partial u}{\partial t}(\cdot,t)\big|_{C^{\alpha}(M,N)}\leq C,
\end{align*}
where the constant \(C\) depends on \(M,N,Z,\nabla Z,V,|d\phi_0|\) and \(T_{max}\).
\end{Lem}
\begin{proof}
Interpreting \eqref{flow-rq} as an elliptic equation we may apply elliptic Schauder theory (\cite{MR925006}, p. 79)
and get that \(u\in C^{1+\alpha}(M,N)\) by the bounds \eqref{bound-dphi} and \eqref{bound-phit}.
Making use of the regularity gained from elliptic Schauder theory we interpret \eqref{flow-rq}
as a parabolic equation. The result then follows by application of parabolic Schauder theory (\cite{MR925006}, p. 79).
\end{proof}

\begin{Lem}
\label{lem-uniqueness}
Let \(u,v\colon M\times [0,T_{max})\to\R^q\) be smooth solutions of \eqref{flow-rq}.
If \(\frac{1}{2}|Z|_{L^\infty}^2\leq\kappa^N\) then the following inequality holds for all \(t\in[0,T_{max})\)
\begin{align*}
|u_t-v_t|^2\leq |u_0-v_0|^2e^{Ct}
\end{align*}
for some positive constant \(C\).
In particular, if \(u_0=v_0\) then \(u_t=v_t\) for all \(t\in[0,T_{max})\).
\end{Lem}

\begin{proof}
In the following \(C\) will denote a universal constant that may change from line to line.
We set \(h:=\frac{1}{2}|u-v|^2\). 
By projecting to a tubular neighborhood \(\sff(du,du), Z(du(e_1)\wedge du(e_2))\) and \(\nabla V(u)\) can be thought of as vector-valued functions in \(\R^q\),
for more details see \cite{MR1896863}, p. 132 and also \cite{MR1800592}, p. 557.
Exploiting this fact a direct computation yields
\begin{align*}
\frac{\partial h}{\partial t}=&\Delta h-|dh|^2-\langle\sff_u(du,du)-\sff_v(dv,dv),h\rangle \\
&-\langle Z_u(du(e_1)\wedge du(e_2))-Z_v(dv(e_1)\wedge dv(e_2)),h\rangle-R\langle\nabla V(u)-\nabla V(v),h\rangle.
\end{align*}
Rewriting
\begin{align*}
\sff_u(du,du)-\sff_v(dv,dv)=(\sff_u-\sff_v)(du,du)+\sff_v(du-dv,du)+\sff_v(dv,du-dv)
\end{align*}
and applying the bounds \eqref{bound-dphi}, \eqref{bound-phit} we find
\begin{align*}
|\langle \sff_u(du,du)-\sff_v(dv,dv),u-v\rangle|&\leq C(|u-v|^2+|du-dv||u-v|), \\
|\langle Z_u(du(e_1)\wedge du(e_2))-Z_v(dv(e_1)\wedge dv(e_2)),u-v\rangle|&\leq C(|u-v|^2+|du-dv||u-v|)
\end{align*}
and the second inequality follows similarly.
Moreover, by assumption \(V\in C^{2,1}(N)\) and thus
\[
|\langle\nabla V(u)-\nabla V(v), u-v\rangle|\leq Ch.
\]
Consequently, by Young's inequality we get
\[
\frac{\partial h}{\partial t}\leq\Delta h+Ch
\]
for some positive constant \(C\). The result then follows by application of the maximum principle.
\end{proof}

\begin{Prop}
Let \(\phi_t\colon M\times [0,T_{max})\to N\) be a smooth solution of \eqref{flow-mfd}.
If \(\frac{1}{2}|Z|_{L^\infty}^2\leq\kappa^N\) then there exists a unique smooth solution for all \(t\in[0,\infty)\).
\end{Prop}
\begin{proof}
This follows from the continuation principle for parabolic partial differential equations.
Suppose that there would be a maximal time of existence \(T_{fin}\), then using the estimates \eqref{bound-dphi} and \eqref{bound-phit}
one can show that we can continue the solution for some small \(\delta>0\) up to \(T_{fin}+\delta\) yielding a contradiction.
The uniqueness follows from Lemma \ref{lem-uniqueness}.
\end{proof}

To achieve convergence of the evolution equation \eqref{flow-mfd} we will make use of
the following 
\begin{Lem}
\label{maximum-principle-l2}
Assume that \((M,h)\) is a compact Riemannian manifold. If a function \(u(x,t)\geq 0\) satisfies
\[
\frac{\partial u}{\partial t}\leq \Delta u+Cu
\]
and if in addition we have the bound
\[
U(t)=\int_Mu(x,t)dM\leq U_0,
\]
then there exists a uniform bound on 
\[
u(x,t)\leq e^CKU_0 
\]
with the constant \(K\) depending on \(M\).
\end{Lem}
\begin{proof}
A proof can for example be found in \cite{MR2744149}, p. 284. For more details on how the 
constant \(C\) in the estimate depends on geometric data, see \cite{MR925006}, Lemma 2.3.1.
\end{proof}

\begin{Prop}[Convergence]
Let \(\phi_t\colon M\times [0,\infty)\to N\) be a smooth solution of \eqref{flow-mfd}.
If \(\Omega\) is exact, \(\frac{1}{2}|Z|_{L^\infty}^2\leq\kappa^N\) and \(|B|_{L^\infty}<1/2\) then the evolution equation \eqref{flow-mfd} 
subconverges in \(C^2(M,N)\) to a harmonic map with scalar and two-form potential \(\phi_\infty\), which is homotopic to \(\phi_0\).
\end{Prop}
\begin{proof}
Since \eqref{flow-mfd} is the negative \(L^2\) gradient flow of the functional \eqref{energy-functional} we obtain
\begin{align*}
\int_M(\frac{1}{2}|d\phi_t|^2+\phi_t^\ast B+RV(\phi_t))dM+\int_0^\infty\int_M\big|\frac{\partial\phi_t}{\partial t}\big|^2dMdt=E(\phi_0).
\end{align*}
Thus, for \(|B|_{L^\infty}<1/2\) we get
\begin{align}
\label{l2-bounds}
\int_M|d\phi_t|^2dM+\int_0^\infty\int_M\big|\frac{\partial\phi_t}{\partial t}\big|^2dMdt\leq C,
\end{align}
where the constant \(C\) depends on \(|B|_{L^\infty},|R|_{L^\infty}|V|_{L^\infty}\) and \(E(\phi_0)\).
Applying Lemma \ref{maximum-principle-l2} to \eqref{estimate-dphi} and \eqref{l2-bounds} we thus obtain a uniform bound on \(|d\phi_t|^2\).
Inserting this bound into \eqref{estimate-phit} we find (for some positive constant \(C\))
\begin{equation*}
\frac{\partial}{\partial t}\frac{1}{2}\big|\frac{\partial\phi_t}{\partial t}\big|^2\leq\Delta\frac{1}{2}\big|\frac{\partial\phi_t}{\partial t}\big|^2
+C\big|\frac{\partial\phi_t}{\partial t}\big|^2.
\end{equation*}
Integrating this equation over \(M\) and \(t\) from \(0\) to \(\infty\) we get 
\begin{equation}
\label{l2-phi-t}
\int_M\big|\frac{\partial\phi_t}{\partial t}\big|^2dM\leq C\int_0^\infty\int_M\big|\frac{\partial\phi_t}{\partial t}\big|^2dMdt
+\int_M\big|\frac{\partial\phi_t}{\partial t}\big|^2\big|_{t=0}dM.
\end{equation}
At this point we may again apply Lemma \ref{maximum-principle-l2} to \eqref{estimate-phit} and \eqref{l2-phi-t}, which gives
\begin{equation*}
\big|\frac{\partial\phi_t}{\partial t}\big|^2_{L^\infty(M\times [0,\infty))}\leq C.
\end{equation*}
Now, by \eqref{l2-bounds}, there exists a sequence \(t_k\) such that
\begin{equation*}
\int_M\big|\frac{\partial\phi_t}{\partial t}(\cdot,t_k)\big|^2dM\to 0
\end{equation*}
as \(k\to\infty\). Moreover, by Lemma \ref{lemma-regularity} we have
\[
\sup_{k\in\N}\big(\big|\frac{\partial\phi_t}{\partial t}(\cdot,t_k)\big|_{C^\alpha(M,N)}+|\phi_{t_k}|_{C^{2+\alpha}(M,N)}\big)\leq C.
\]
Then it follows by the Theorem of Arzela and Ascoli that there exists a convergent subsequence, which we will
also denote by \(t_k\). Thus, the family \(\phi_{t_k}\) subconverges in \(C^2(M,N)\) to a map \(\phi_\infty\). Since \(\phi_t\)
depends smoothly on \(t\), the limit \(\phi_\infty\) is homotopic to \(\phi_0\).
\end{proof}

\begin{Bem}
Note that we only made use of the globally defined energy \eqref{energy-functional} to obtain convergence
of the gradient flow. To achieve long-time existence of \eqref{flow-mfd} we do not require a variational structure.
\end{Bem}

\subsection{Non-compact Target}
In the case that the target manifold \(N\) is complete, but non-compact we have to make additional assumptions to control
the image of \(M\) under the evolution of \(\phi\). However, we can use the potential \(V(\phi)\) to constrain \(\phi_t(M)\) 
to a compact set in \(N\). This is similar to the case of the heat flow for harmonic maps with potential \cite{MR1800592}.
To this end, let \(d_N(y)\) denote the Riemannian distance in \(N\) from some fixed point \(y_0\).

Finally, we will prove the following:
\begin{Satz}
\label{theorem-non-compact}
Let \((M,h)\) be a closed Riemann surface and \((N,g)\) a complete, oriented Riemannian manifold
with negative sectional curvature. Moreover, suppose that \(\Omega\) is exact, \(|B|_{L^\infty}<1/2\), \(V\in C^{2,1}(N)\)
and \(\frac{1}{2}|Z|_{L^\infty}^2\leq\kappa_N\), where \(\kappa_N\) denotes and upper bound on the sectional curvature on \(N\).  
In addition, assume that the potential \(V(\phi)\) satisfies
\begin{equation*}
-R\hess V(y)\leq-\frac{C}{1+d_N(y)} \qquad\text{ or }\qquad |R\hess V|_{L^\infty}\leq\frac{\lambda_1(M)}{2},
\end{equation*}
where \(\lambda_1(M)\) denotes the first eigenvalue of the Laplacian on \(M\).
Let \(\phi_t\colon M\times [0,\infty)\to N\) be a smooth solution of \eqref{flow-mfd}.
Then \(\phi_t\) subconverges in \(C^2(M,N)\) to a harmonic map with scalar and two-form potential,
which is homotopic to \(\phi_0\). Moreover, \(\phi_\infty\) is minimizing the energy in its homotopy class.
\end{Satz}

First of all, we make the following observation:
\begin{Lem}
If the three-form \(\Omega\) is exact, then \(Z \in \Gamma (\Hom(\Lambda^2T^\ast N,TN))\) defined via \eqref{def-Z}
is parallel, that is \(\nabla Z=0\).
\end{Lem}
\begin{proof}
We fix a point \(p\in N\) and we extend any tangent vectors \(\xi_1,\xi_2,\eta\in T_pN\) to vector fields being
\emph{synchronous} in the point \(p\), meaning that
\[
\nabla_V\xi_1=\nabla_V\xi_2=\nabla_V\eta=0
\]
for all tangent vectors \(V\in T_pN\). Then using that \(\Omega=dB\) is exact, we find from \eqref{def-Z}
\[
dB(\xi_1,\xi_2,\eta)=\langle Z(\xi_1\wedge\xi_2),\eta\rangle.
\]
By differentiation we then obtain
\[
0=\langle\nabla Z(\xi_1\wedge\xi_2),\eta\rangle,
\]
which proves the claim.
\end{proof}

As a next step we show how we can use the Hessian of the scalar potential \(V(\phi)\) to
constrain \(\phi_t(M)\) to a compact set. This idea has already been applied in the study of the 
heat flow for harmonic maps with potential, see \cite{MR1800592}, Proposition 2.
\begin{Lem}
Let \((M,h)\) be a closed Riemann surface and \((N,g)\) a complete Riemannian manifold.
Suppose that \(\Omega\) is exact and that \(\frac{1}{2}|Z|_{L^\infty}^2\leq\kappa_N\).
Moreover, assume that the potential \(V(\phi)\) satisfies
\begin{equation}
\label{lemma-assumption-potential-1}
-R\hess V(y)\leq-\frac{C}{1+d_N(y)}.
\end{equation}
Then there exists a compact set \(K\subset N\) such that \(\phi_t(M)\subset K\) as long as the solution of \eqref{flow-mfd} exists.
\end{Lem}
\begin{proof}
Let \(d_N\) be the distance function from some point \(y_0\) in \(N\). We set
\[
f(t):=\sup_{x\in M}\big|\frac{\partial\phi_t}{\partial t}\big|.
\]
It follows directly that
\begin{equation*}
d_N(\phi_t)\leq C+\int_0^tf(s)ds.
\end{equation*}
Making use of the assumptions on the Hessian of the potential \eqref{lemma-assumption-potential-1}, the homomorphism \(Z\) and the sectional curvature \(K^N\),
we obtain from \eqref{estimate-phit} by an argument similar to the compact case that
\[
\frac{\partial}{\partial t}\big|\frac{\partial\phi_t}{\partial t}\big|^2\leq\Delta\big|\frac{\partial\phi_t}{\partial t}\big|^2
-\frac{C}{1+\int_0^tf(s)ds}\big|\frac{\partial\phi_t}{\partial t}\big|^2.
\]
Then, by the maximum principle, we find
\begin{equation*}
f(t)\leq f(0)\exp\big(-\frac{C}{2}\int_0^t\frac{1}{1+\int_0^\tau f(s)ds}d\tau\big).
\end{equation*}
We can deduce that
\begin{equation*}
\int_0^tf(s)ds\leq C,
\end{equation*}
which gives
\begin{equation*}
f(t)\leq f(0)e^{-Ct}.
\end{equation*}
Hence, we find that \(d_N(y)\leq C\). This proves that there exists a compact set \(K\subset N\)
such that \(\phi_t(M)\subset K\) for all \(t\in [0,T_{max})\).
\end{proof}

As noted in \cite{MR1979036}, Proposition 2.4, one can also constrain \(\phi_t(M)\) to a compact set
if the maximum of the Hessian of the potential \(V(\phi)\) is bounded by the first eigenvalue of
the Laplacian on \(M\). This idea can also be applied here:

\begin{Lem}
Let \((M,h)\) be a closed Riemann surface and \((N,g)\) a complete, oriented Riemannian manifold.
Suppose that \(\Omega\) is exact and that \(\frac{1}{2}|Z|_{L^\infty}^2\leq\kappa_N\).
Moreover, assume that the potential \(V(\phi)\) satisfies
\begin{equation}
\label{lemma-assumption-potential-2}
|R\hess V|_{L^\infty}\leq\frac{\lambda_1(M)}{2},
\end{equation}
where \(\lambda_1(M)\) denotes the first eigenvalue of the Laplacian on \(M\).
Then there exists a compact set \(K\subset N\) such that \(\phi_t(M)\subset K\) as long as the solution of \eqref{flow-mfd} exists.
\end{Lem}
\begin{proof}
Using \eqref{bochner-phit} we derive the following inequality
\begin{align*}
\frac{\partial}{\partial t}\frac{1}{2}\int_M\big|\frac{\partial\phi_t}{\partial t}\big|^2dM\leq &\int_M\big(-\big|\nabla\frac{\partial\phi_t}{\partial t}\big|^2
-\kappa_N|d\phi_t|^2\big|\frac{\partial\phi_t}{\partial t}\big|^2+
|Z|_{L^\infty}\big|\nabla\frac{\partial\phi_t}{\partial t}\big||d\phi_t|\big|\frac{\partial\phi_t}{\partial t}\big| \\
\nonumber &-R\hess V(\frac{\partial\phi_t}{\partial t},\frac{\partial\phi_t}{\partial t})\big)dM \\
\nonumber \leq&\int_M\big(-\frac{1}{2}|\nabla\frac{\partial\phi_t}{\partial t}\big|^2+|R\hess V|_{L^\infty}\big|\frac{\partial\phi_t}{\partial t}\big|^2\big)dM,
\end{align*}
where we used Young's inequality and the assumptions on the homomorphism \(Z\) and the sectional curvature \(K^N\) in the last step.
By the Kato inequality \(\big|\nabla f\big|^2\geq \big|\nabla|f|\big|^2\) for a function \(f\colon M\to\R\) and the Poincaré inequality on \(M\)
we obtain
\begin{equation*}
\int_M|\nabla f|^2dM\geq\lambda_1(M)\int_M f^2dM.
\end{equation*}
Applying this inequality  we find
\[
\frac{\partial}{\partial t}\frac{1}{2}\int_M\big|\frac{\partial\phi_t}{\partial t}\big|^2dM\leq(|R\hess V|_{L^\infty}-\frac{\lambda_1(M)}{2})\int_M\big|\frac{\partial\phi_t}{\partial t}\big|^2dM
\]
yielding
\begin{equation*}
\int_M\big|\frac{\partial\phi_t}{\partial t}\big|^2dM\leq C_0e^{(2|R\hess V|_{L^\infty}-\lambda_1(M))t}.
\end{equation*}
Using the assumption on \(|R\hess V|_{L^\infty}\) and applying Lemma \ref{maximum-principle-l2} to \eqref{bochner-phit} 
we obtain a uniform pointwise bound on \(|\frac{\partial\phi_t}{\partial t}|\). By the same argument as in the previous Lemma
this yields the claim.
\end{proof}

By making use of the previous Lemmata we thus find

\begin{Lem}
Let \((M,h)\) be a closed Riemann surface and \((N,g)\) a complete, oriented Riemannian manifold.
Moreover, suppose that \(\Omega\) is exact, \(|B|_{L^\infty}<1/2\), \(V\in C^{2,1}(N)\)
and \(\frac{1}{2}|Z|_{L^\infty}^2\leq\kappa_N\).
In addition, assume that the potential \(V(\phi)\) satisfies
\begin{equation*}
-R\hess V(y)\leq-\frac{C}{1+d_N(y)} \qquad\text{ or }\qquad |R\hess V|_{L^\infty}\leq\frac{\lambda_1(M)}{2},
\end{equation*}
where \(\lambda_1(M)\) denotes the first eigenvalue of the Laplacian on \(M\).
Let \(\phi_t\colon M\times [0,\infty)\to N\) be a smooth solution of \eqref{flow-mfd}.
Then \(\phi_t\) converges in \(C^2(M,N)\) to a harmonic map with scalar and two-form potential,
which is homotopic to \(\phi_0\). 
\end{Lem}
\begin{proof}
By the previous Lemmata we know that \(\phi_t(M)\) stays inside a compact set \(K\). Thus, by Theorem \ref{theorem-compact} there exists
a sequence \(t_k\) such that \(\phi_{t_k}\) converges to a harmonic map with scalar and two-form potential.
Moreover, we have
\[
d_N(\phi_t,\phi_\infty)\leq d_N(\phi_t,\phi_{t_k})+d_N(\phi_{t_k},\phi_\infty).
\]
We know that
\begin{equation*}
d_N(\phi_t,\phi_{t_k})\leq\int_{t_k}^t\big|\frac{\partial\phi_s}{\partial s}\big|ds\leq C\int_{t_k}^te^{-Cs}ds\to 0
\end{equation*}
as both \(k,t\to\infty\), which proves the claim.
\end{proof}

\begin{Bem}
In the case of the heat flow for harmonic maps with potential to a non-compact target 
with a concave potential the limit \(\phi_\infty\) is trivial, see \cite{MR1800592}, p. 564.
This statement relies on the Bochner formula \eqref{bochner-vphi}, due to the presence of the two-form potential
we cannot draw the same conclusion here.
\end{Bem}

\begin{Bem}
If we would only require an upper bound on \(R\hess V\) in \(N\) then this would be enough to establish
long-time existence of \eqref{flow-mfd}. Consequently, to achieve convergence of \eqref{flow-mfd} we need
the potential \(RV(\phi)\) to constrain \(\phi_t(M)\) to a compact set.
\end{Bem}

\begin{Bem}
Using the extrinsic version of the evolution equation \eqref{flow-rq} we can apply the maximum principle to bound the image of \(\phi_t(M)\).
This idea was already used for related geometric flows: A similar criterion as below for the heat flow of harmonic maps with
potential is given in \cite{MR1800592}, Proposition 2. For the heat flow of the prescribed mean curvature equation in \(\R^3\) 
the same idea is used in Proposition 3.2 in \cite{MR1125012}.

By a direct calculation using \eqref{flow-rq} we obtain
\begin{align*}
\frac{\partial}{\partial t}|u|^2=&\Delta|u|^2-2|du|^2+\langle u,\sff(du,du)+Z(\uarg)\rangle-R\langle u,\nabla V(u)\rangle \\
\leq &\Delta |u|^2-2|du|^2(1-\frac{|u|}{2}(|\sff|_{L^{\infty}}+|Z|_{L^\infty}))-R\langle u,\nabla V(u)\rangle.
\end{align*}
Hence, if we can guarantee that
\[
\frac{|u|}{2}(|\sff|_{L^{\infty}}+|Z|_{L^\infty})\leq 1, \qquad R\langle u,\nabla V(u)\rangle>0
\]
at \(t=0\), then by the maximum principle we get a bound on \(|u|^2\) as long as the solution exists.
However, the second condition on the potential \(V(u)\) is hard to ensure.
\end{Bem}

\subsection{Minimizing the energy}
In this section we briefly discuss if the limiting map constructed in Theorem \ref{theorem-non-compact}
is minimizing energy in its homotopy class.
\begin{Lem}
Under the assumptions of Theorem \ref{theorem-non-compact} the limiting map \(\phi_\infty\) is minimizing energy in its homotopy class.
\end{Lem}
\begin{proof}
We use the formula for the second variation \eqref{second-variation} to
show that the limit \(\phi_\infty\) obtained in Theorem \ref{theorem-non-compact} is minimizing the energy in its homotopy class.
Let \(\phi_1,\phi_2\colon M\to N\) be two smooth maps and let \(\Phi\colon M\times [0,1]\) be 
a geodesic homotopy between \(\phi_1\) and \(\phi_2\), that is \(\Phi(\cdot,0)=\phi_1\) and
\(\Phi(\cdot,1)=\phi_2\). Moreover, assume that for any \(x\in M\) the map \(\Phi(x,\cdot)\)
is a geodesic, when \(\phi_1\) is a harmonic map with scalar and two-form potential. 
Using the formula for the second variation \eqref{second-variation}
we find the following inequality
\begin{align*}
E(\phi_1)-E(\phi_2)=&\int_0^1d\sigma\int_0^\sigma
\big(\big|\nabla\frac{\partial\Phi}{\partial s}\big|^2+\kappa_N|d\Phi|^2\big|\frac{\partial\Phi}{\partial s}\big|^2
-|Z|_{L^\infty}|d\Phi|\big|\nabla\frac{\partial\Phi}{\partial s}\big|\big|\frac{\partial\Phi}{\partial s}\big|\\
&-R\hess V(\frac{\partial\Phi}{\partial s},\frac{\partial\Phi}{\partial s})
\big)ds \\
&\geq\int_0^1d\sigma\int_0^\sigma
\big(\frac{1}{2}\big|\nabla\frac{\partial\Phi}{\partial s}\big|^2+(\kappa^N-\frac{1}{2}|Z|^2_{L^\infty})|d\Phi|^2\big|\frac{\partial\Phi}{\partial s}\big|^2 \\
&-R\hess V(\frac{\partial\Phi}{\partial s},\frac{\partial\Phi}{\partial s})\big)ds \\
>&0.
\end{align*}
Thus, under the assumptions of the Lemma \(\phi_1\) is an energy minimizer in its homotopy class.
\end{proof}

\begin{Bem}
Let \(\phi_t\colon M\times [0,\infty)\to N\) be a smooth solution of \eqref{flow-mfd}.
Under the assumptions of Theorem \ref{theorem-non-compact} the energy \(E(\phi(t))\) is a convex function of \(t\),
which follows by a direct calculation.
\end{Bem}

Exploiting the convexity of the energy \(E(\phi(t))\) with respect to \(t\) we can prove the following
uniqueness Theorem, which is very similar to the case of Hartman's theorem for harmonic maps \cite{MR0214004}.

\begin{Prop}
Under the assumptions of Theorem \ref{theorem-non-compact} the limit \(\phi_\infty\)
is independent of the chosen subsequence \(t_k\).
\end{Prop}
\begin{proof}
The proof is the same as in the case of harmonic maps and we omit it here.
\end{proof}

\bibliographystyle{plain}
\bibliography{mybib}
\end{document}